\documentclass[11pt]{article}
\usepackage[margin=2.5cm]{geometry}
\usepackage[affil-it]{authblk}
\usepackage[utf8]{inputenc}
\usepackage{amsmath,amsthm,amssymb}
\usepackage{float}
\usepackage{tikz}
\usepackage{keyval}
\usepackage{ifthen}
\usepackage{booktabs}
\usepackage{epsf}
\usepackage{psfrag}
\usepackage{wrapfig}
\usepackage{color,rotating,psfrag,amsfonts}
\usepackage{todonotes}
\usepackage[dvips]{epsfig}
\usepackage{subfigure}
\usepackage{adjustbox}

\usepackage{tikz-network}

\usepackage{hyperref} % for url

\usepackage{amssymb,amsmath, amsthm}

\usepackage[title]{appendix}

\usepackage{graphicx}
\usepackage{pdfpages}
\usepackage{hyperref}
\usepackage{mathtools}
\usepackage[capitalise]{cleveref}
\usepackage{nicefrac}
\usepackage{stmaryrd}
\usepackage{authblk}
\usepackage{subcaption}
\usepackage[pagewise]{lineno}

\usepackage{changes}

\title{Multi-objective optimization for multi-agent injection strategies in subsurface CO$_2$ storage}

\begin{document}
\date{}
%\author[$\star, 1$]{Per Pettersson}
%\author[$\dagger$]{Sebastian Krumscheid}
%\author[$\star,\ddagger$]{Sarah Gasda}
% \author[$\star$]{Alfred Schmidt}

%\affil[$\star$]{NORCE Norwegian Research Centre, N5838 Bergen, Norway}
%\affil[$\dagger$]{Karlsruhe Institute of Technology (KIT), Kaiserstraße 12, 
%Karlsruhe, 76131, Germany}

%\affil[$\ddagger$]{University of Bergen, 5020 Bergen, Norway}

\author[1]{Per Pettersson\thanks{Corresponding author: per.pettersson@norceresearch.no}}

\author[2]{Sebastian Krumscheid}
\author[1,3]{Sarah Gasda}

\affil[1]{NORCE Norwegian Research Centre, N-5838 Bergen, Norway}
\affil[2]{Karlsruhe Institute of Technology (KIT), 76131 Karlsruhe, Germany}
\affil[3]{University of Bergen, N-5020 Bergen, Norway}

%\author[inst1,inst3]

%\affiliation[inst1]{organization={NORCE Norwegian Research Centre},%Department and Organization
%            addressline={PO Box 22, Nygardstangen}, 
%            city={Bergen},
%            postcode={5838}, 
            %state={State One},
%            country={Norway}}

%\affiliation[inst2]{organization={Karlsruhe Institute of Technology (KIT)},%Department and Organization
%            addressline={Kaiserstraße 12}, 
%            city={Karlsruhe},
%            postcode={76131}, 
            %state={State Two},
%            country={Germany}}   
            
%\affiliation[inst3]{organization={University of Bergen},%Department and Organization
%            addressline={PO Box 7803}, 
%            city={Bergen},
%            postcode={5020}, 
         %   state={State Two},
%            country={Norway}}

\maketitle

\begin{abstract}
%Efficient use of available reservoir resources is key in large-scale subsurface CO$_2$ storage. The optimal injection strategy for a given set of injection wells is affected not only by the supply of CO$_2$ and the geology, but also by activities performed by external parties leading to pressure buildup, e.g., competing companies that inject CO$_2$ within the same or a geologically connected formation. Assuming that all agents performing CO$_2$ injection aim to maximize their own total storage amounts, and are free to form binding agreements with other agents, it is an open question what strategies they should choose to achieve the optimal outcome.
%In this work, we 
We propose a  novel framework for optimizing injection strategies in large-scale CO$_2$ storage  combining multi-agent models with multi-objective optimization, and reservoir simulation. We investigate whether agents, i.e., well operators, should form coalitions for collaboration to maximize the outcome of their storage activities. 
In multi-agent systems, it is typically assumed that the optimal strategy for any given coalition structure is already known, and it remains to identify which coalition structure is optimal according to some predefined criterion. In the setting considered here, for any given coalition structure, the optimal CO$_2$ injection strategy is not a priori known, and needs to be found by a combination of reservoir simulation and optimization. Since all coalitions seek to maximize their own storage volumes, the optimal strategy for every coalition structure is the solution to a multi-objective optimization problem. 

The multi-objective optimization problems all come with the numerical challenges of repeated evaluations of complex-physics models. We use versatile evolutionary algorithms to solve the sets of multi-objective optimization problems, where the solution is a set of values, e.g., a Pareto front. The Pareto fronts are first computed using the so-called weighted sum method that transforms the multi-objective optimization problem into a set of single-objective optimization problems. Results based on two different Pareto front selection criteria are presented. Then a truly multi-objective optimization method is employed to obtain the Pareto fronts, and we investigate its performance compared to the previous weighted sum method. 

We demonstrate the proposed framework on the Bjarmeland formation, a pressure-limited prospective storage site in the Barents Sea. The problem is constrained by the maximum sustainable pressure buildup and an assumed supply of CO$_2$ that can vary over time. In addition to identifying the optimal coalitions, the methodology also shows how distinct suboptimal coalitions perform in comparison to the optimum.

\end{abstract}



%\begin{keyword}
%keyword one \sep keyword two
%\end{keyword}

%\end{frontmatter}

%%%%%%%%%%%%%%

\section{Introduction}
%%%%
Large-scale CO$_2$ storage is a key technology to reach the goals described by the Intergovernmental Panel on Climate Change (IPCC) to significantly reduce the projected human emissions of greenhouse gases to the atmosphere. Limiting global warming to at most 1.5 $^{\circ}$C requires annual injection of 3-10 Gt within the next few decades~\cite{IPCC_23}. With maturing technology for achieving CO$_2$ injection at such large scales, multi-site and basin-scale utilization of subsurface CO$_2$ storage resources become feasible, while socio-economic challenges remain and need further investigation~\cite{Krevor_etal_23}. In addition to the risk of directly compromising reservoir integrity by pressure buildup and other forms of potential leakage, CO$_2$ storage operations are likely to be affected by (and themselves affect) natural gas storage and hydrocarbon extraction in hydraulically-connected reservoir sites. Efficient utilization of large-scale CO$_2$ sites therefore calls for methods to optimize storage where both physical and socio-economic effects are taken into account. In particular, multi-agent systems models for CO$_2$ injection, to be introduced in this paper, are relevant to account for the effect of independent agents with conflicting goals and make sure that all agents attain their respective goals in terms of value and risk avoidance.

Basin-scale CO$_2$ storage estimation via simulation is challenging for a number of reasons, including vast physical and temporal scales, and complex physics on different scales. 
Pressure buildup often imposes severe constraints on the amounts of CO$_2$ that can be safely injected without jeopardizing the integrity of the reservoir~\cite{Rutqvist_etal_07}.
For modeling of basin-scale CO$_2$ storage with multiple injection wells, pore-volume estimates of total storage capacity are not accurate enough and need to be complemented by simulation of pressure buildup and associated consequences on caprock integrity and potential leakage to groundwater resources~\cite{Birkholzer_Zhou_09}. Analytical pressure buildup models by means of superpositions of analytical or semi-analytical solutions cannot capture structural variability, as demonstrated in~\cite{Huang_etal_14}, where pressure buildup models of different complexity were compared. Rare-event simulation of extreme event CO$_2$ capacity estimates on the regional scale was performed in~\cite{Pettersson_etal_22}, demonstrating the feasibility of accurate uncertainty quantification methods at these scales.
A framework for regional-scale CO$_2$ injection in multiple wells combining static approximation with dynamical simulation was demonstrated in~\cite{Tveit_etal_24}.

%{\color{blue}Reduced-physics vs surrogate models:}\\
Optimization of CO$_2$ storage requires a large number of evaluations of physical models and becomes infeasible unless models with limited computational cost and acceptable accuracy are employed.
To remedy the problem of prohibitive computational cost of repeated model evaluations in optimization, the discretized partial differential equations model can be partly replaced by a surrogate model, e.g., as demonstrated for hysteretic trapping to optimize well placement~\cite{Babaei_etal_15}, and using hierarchical interpolation for microbially induced calcite precipitation in~\cite{Tveit_etal_20}. Perhaps more common in the field of CO$_2$ storage simulation is to replace a full-scale complex physical model with a simplified-physics model, where the problem has been carefully analyzed so that less important physical mechanisms can be ignored in favor of increased resolution of more significant physical mechanisms. In particular, vertically-integrated models have gained considerable interest and successfully been demonstrated to serve as realistic large-scale reservoir models, c.f.~\cite{Gasda_etal_09, Nilsen_etal_11, Guo_etal_14}. %{\color{blue}[@Sarah: any preferred VE references here? General CO2 VE models can go here.]}

%{\color{blue}Opt of CO2 with simplified models below:}\\
Optimization methods can be divided into global optimization methods that are general but often expensive, versus local, gradient-based methods that require knowledge or approximation of the gradients of the objective and constraint functions.
Within the context of CO$_2$ storage, gradient-based methods have previously been employed for well rate optimization~\cite{Allen_etal_17b}, while global optimization methods with their reduced risk of getting stuck in a local optimum have been used for well location optimization and optimization of multiple properties where gradients are not easily available~\cite{Cameron_Durlofsky_12, Cihan_etal_15, Musayev_etal_23}.

%{\color{blue}[Reiterate the imprortance of pressure constraint.]}
For pressure-limited prospective CO$_2$ storage sites, pressure will act as a constraint that needs to be included in the optimization methods for fully utilizing the storage potential.
Optimization with an analytical single-phase model was used for pressure management by optimization of fluid (e.g., brine) extraction wells to mitigate pressure buildup in~\cite{Birkholzer_etal_12}. In addition to reduction of computational time, a simplified physical model can also simplify the optimization problem itself, as exemplified in~\cite{Santibanez_etal_19}, where simplification to linear objective functions and constraints resulted in a linear programming problem solved with the Simplex algorithm.

The work described so far relies on the tacit assumption that all operations affecting a storage site should be simultaneously optimized. This is realistic if a single agent, e.g., a big company, controls all operations on hydraulically connected storage sites. In contrast, if independent agents perform operations that affect each other, e.g., by means of pressure buildup limiting injection,  they are of course very unlikely to willfully underperform in their own activities to the benefit of competing companies. In this paper, we propose a multi-agent model for CO$_2$ storage to account for distinct agents that have their own goals, and can form coalitions for collaboration under binding agreements for the mutual benefits of the members of the coalition. The agents can act independently but will still be affected by other agents in their injection operations. Moreover, the agents are assumed to be rational in pursuing their goals, and their actions are determined by the simultaneous optimization of conflicting objectives. This requires a multi-objective optimization of CO$_2$ site operations and comes with the challenges of standard optimization, and some additional challenges to be addressed. Work on multi-objective optimization for CO$_2$ storage has so far been limited and within the setting of a single agent that wishes to fulfill more than one objective. Storage volume of CO$_2$ was considered one of the objectives in a multi-objective optimization workflow for reservoir simulation, where an artificial neural network was used as a surrogate model to replace expensive compositional reservoir models~\cite{You_etal_20}.
A non-dominated sorting genetic algorithm was used with a heterogeneous aquifer model to simultaneously optimize the two objectives cost, defined by means of injection pressure, and CO$_2$ containment efficiency~\cite{Park_etal_21}.

In this paper, we propose a multi-agent model for multi-criteria decision-making for CO$_2$ storage, to be described in the next section. This model leads to a set of potential constellations for collaboration, referred to as coalition structures, and the performance in terms of successfully injecting CO$_2$ is determined by the solutions to a multi-objective optimization problem for every potential coalition structure.
These multi-objective optimization problems all come with the challenges of repeated evaluations of complex-physics models. We use evolutionary algorithms to solve the sets of multi-objective optimization problems, and visualize Pareto fronts depicting these solutions for the Bjarmeland formation, a pressure-limited prospective storage site in the Barents Sea.
The Pareto fronts are computed using the so-called weighted sum method that transforms the multi-objective optimization problem into a set of single-objective optimization problems. Results based on two different Pareto front selection criteria are presented. Then a truly multi-objective optimization method is employed to obtain the Pareto fronts, and we investigate its performance compared to the previous weighted sum method. Finally, we provide a discussion of the results and conclusions.

\section{Multi-agent Model for CO$_2$ Injection}

Consider a situation with a large-scale CO$_2$ storage site with injection and possibly extraction wells that are operated by different agents with individual and often conflicting objectives. In particular, pressure buildup resulting from the injection operations of a given agent can have adverse effects on the realizable injection rates of the other agents. The agents are typically commercial companies that strive to maximize their own total CO$_2$ injections, and are free to form agreements for collaboration with other agents, but are equally free not to do so. The union of operators of two or more wells that join forces to optimize the total performance of the wells, rather than individually trying to optimize for a single well, will be referred to as a \textit{coalition}. For the modeling framework to be proposed here, it does not matter whether two wells in a coalition are otherwise operated by the same agent, or by two distinct agents. 
Individual wells that are operated without cooperation with other wells are called singleton coalitions, and the set of all coalitions is denoted \textit{coalition structure}. Note that a coalition structure in the current setting is a partition of all wells into disjoint coalitions, so every well belongs to exactly one coalition. The number of possible coalitions for a set of $a$ agents is $2^{a}-1$, and the number of possible coalitions structures, known as the Bell number, grows even faster with $a$. For instance, for $a=3,4,5$, there are respectively 5, 15, and 52 coalition structures, and already for $a=10$ there are as many as 115975 coalition structures. 
As an example, with three wells (W1, W2, W3), there are the five coalition structures: (i) \{\{W1,W2,W3\}\} (grand coalition), (ii) \{\{W1,W2\}, \{W3\}\}, (iii) \{\{W1,W3\}, \{W2\}\}, (iv) \{\{W1\}, \{W2,W3\}\}, and (v) \{\{W1\}, \{W2\}, \{W3\}\} (singleton coalitions). If one wants to investigate properties pertaining to a certain coalition structure that can only be obtained at some non-negligible computational expense, i.e., the optimal injection schedule for that coalition structure, it is clear that the number of agents needs to remain small. However, even for a large-scale CO$_2$ storage site, the number of distinct agents is not likely to be very large, and there may be certain coalition structures that are a priori not feasible and hence not relevant for further investigation. Hence, the complexity growth may in some practical situations be less severe than what the Bell number indicates.

A coalition structure over the set of all $a$ agents describes a tentative structure for collaboration, where it is assumed that the members of the coalitions seek to optimize the total value of the coalition. Throughout the remainder of the paper, the \emph{value} of a coalition is defined as the total amount of CO$_2$ the coalition can safely inject, subject to physical constraints. Other definitions of value, e.g., risk avoidance measures or measures combining different desired outcomes, can be used as long as they can be inferred from a given model. One should note that even in cases where the individual contributions from agents to a coalition are clear, e.g., how much CO$_2$ has been injected in any individual well, the distribution of value between the coalition members need not necessarily correspond to the individual contributions. For instance, if a given agent has injected less than physically possible to maximize the total amount of the coalition, that agent should probably be assigned a value larger than its actual contribution as an incentive to not leave the coalition.

Multi-criteria decision-making provides a systematic framework for choosing the most suitable candidate among a potentially large set of possible solutions, each one to a varying degree satisfying a number of desired criteria to be fulfilled~\cite{Triantaphyllou_00}. In essence, the framework provides decision-making algorithms that are assumed to rely on inputs and outputs from a decision maker with a set of preferences and goals to be satisfied. Depending on how the decision-maker articulates preferences, different methods are available for forming a decision. In \textit{a priori} preference articulation, the decision maker has stated preferences from the start, while in progressive preference articulation, the decision maker sequentially provides preference input to the decision-making algorithm. In the current work, we consider \textit{a posteriori} decision-making, where we seek to provide a wide range of solution outputs for the decision maker to choose from after the decision-making algorithm has finished. Advantages of this approach include flexibility in the sense that the decision maker can investigate multiple decision paths, as well as making more informed decisions by incorporating the knowledge gained from the output of the algorithm. Clearly, the preferred coalition structure depends on the relative preferences of the decision maker. One should note that in the context of independent agents, the concept of decision maker can be ambiguous. An actual agreement, or a series of independent practical considerations, can in practice lead to the selection being made, and hence define the decision-making process.

In order to determine an ideal coalition structure among all possible coalition structures, i.e., one that in some sense leads to optimal outcomes for all agents so that they do not choose to break the coalition structure and form a more profitable coalition, the values of all coalitions and coalition structures need to be determined. While the assignments of individual values are commonly assumed to be achievable at unit computational cost~\cite{Rahwan_etal_15}, in the current paper numerical optimization with multiple objectives need to be performed. For every possible coalition structure, a multi-objective optimization problem needs to be solved. Next, we briefly describe the framework for multi-objective optimization.

\section{Multi-objective Optimization of CO$_2$ Injection Rates}
%{\color{blue}[Check whether 'objective' is clearly explained/defined in our context!]}
For any coalition structure of size $2 \leq m \leq a$, there are
 $m$ possibly conflicting coalition objectives $F_1(q), F_2(q),\dots, F_{m}(q)$ to be optimized. In this work, every objective corresponds to the value (i.e., total amount of CO$_2$) of a given coalition embedded in a coalition structure and it is  a function of up to $n$ design variables $q = (q_1,q_2,\dots, q_n) \in \mathbb{R}^{n}$ which are the free variables with respect to which we optimize the coalition objectives. The design variables will henceforth correspond to the injection schedules of the agents, i.e., the injection rates for any given well and time interval, and the coalition objectives will be linear functions of the design variables.
 The design variables are subject to a nonlinear inequality constraint function $g = (g_1,g_2,\dots, g_k) \in \mathbb{R}^{k}$ such that $g_i(q) \leq 0$, for $i=1,\dots, k$. 
 Here, the truly nonlinear part of the constraint function $g$ represents an upper limit on the reservoir pressure $p_{\text{res}}$ set to 90~\% of the overburden pressure $p_{\text{ob}}$, evaluated pointwise in space,
 \[
 g_i = \frac{p_{ \text{res}}(\mathbf{x}_{i}) }{p_{\text{ob}}(\mathbf{x}_{i}) } - 0.9, \quad \mbox{for } i \mbox{ so that } \mathbf{x}_{i} \in \mathbf{X}_{\text{grid}},
 \]
 where $\mathbf{X}_{\text{grid}}$ denotes the discrete spatial grid of the numerical reservoir model. There may also be linear inequality constraints on the injection rates, corresponding to minimum economically feasible and maximum supply rates.
 The resulting multi-objective optimization (MOO) problem,
\begin{equation}
\label{eq:moo}
\begin{aligned}
\max_{q} \left(
F_1(q), F_2(q),\dots, F_{m}(q)
\right)\\
\mbox{ s.t. } g(q) \leq 0
\end{aligned},
\end{equation}
typically does not have a unique solution since a unique $q$ that simultaneously maximizes all $F_1,\dots, F_m$ usually does not exist. The special case of $m=1$ represents the grand coalition where all agents choose to collaborate and~\eqref{eq:moo} then reduces to a constrained single-objective optimization problem. In the following, we assume that the constraints $g$ are such that a candidate solution always exists, i.e., the feasible design space $\{q | g(q) \leq 0 \}$ is non-empty. The problem is otherwise infeasible. A candidate solution $q^*$ is said to be Pareto optimal (or Pareto efficient) if it is not possible to improve one objective without deteriorating another objective. %{\color{blue}[Could add some more details on Pareto optimality.]}
We seek to identify the set of Pareto optimal solutions, all of which are in some sense candidates to be selected as the best injection schedules to be implemented by the well operators. The points of the Pareto fronts then represents different optimized injection scenarios. However, some of the solutions may not be desirable by all agents simultaneously. In what follows, we briefly describe methods to approximate these Pareto fronts and discuss the selection of solution candidates.

\subsection{Pareto Front Approximation}
\label{sec:pareto_appr}
Methods to solve~\eqref{eq:moo} can be sufficient, necessary, or both, to yield candidate solutions that are Pareto optimal. A popular and particularly simple method that always yields Pareto optimal solutions (i.e., sufficient condition for Pareto optimality) is the weighted sum method (WSM)~\cite{Zadeh_63}, i.e., maximization of the weighted sum of objective functions,
\begin{equation}
\label{eq:WSM}
F_{\text{WSM}} = \sum_{j=1}^{m} w_j F_j(q),
\end{equation}
where the weights satisfy $\sum_{j=1}^{m} w_j=1$,  $w_j > 0$, and assuming the same constraint function $g$ as in~\eqref{eq:moo}. By varying the weights, a Pareto front can be approximated. The WSM is attractive since it is easy to use and the MOO problem~\eqref{eq:moo} is transformed to a set of constrained single-objective optimization (SOO) problems that are typically easier to solve than the full MOO. In particular, if one does not need a well-resolved Pareto front, only a single or small number of different weights can be considered to limit the total amount of model evaluations. It is however not straightforward to find a direct correspondence between WSM weights and Pareto front location, as investigated and discussed in detail in~\cite{Marler_Arora_10}.

A broad range of evolutionary algorithms have been developed to target multi-objective optimization problems and are widely used in many applications~\cite{Zhou_etal_11}. We use evolutionary algorithms both to solve the set of SOO resulting from applying the WSM, and, as a comparison, to directly compute the Pareto fronts corresponding to~\eqref{eq:moo}. As pointed out in the Introduction, the algorithms in this family are robust, scale well with the number of design variables and are relatively straightforward to extend to different situations not only including CO$_2$ injection rates.
For the SOO problems, we employ the competitive swarm method (CSM). It is a variant of particle swarm methods with enhanced population diversity where the population is divided into a superior half that is transferred directly to the next generation, and an inferior half that is updated based on the superior individuals~\cite{Cheng_Jin_14}. CSO is suitable for high-dimensional problems, and there are versions for both single-objective and multi-objective optimization~\cite{Chauhan_etal_24}.
To solve the MOO~\eqref{eq:moo} directly with a constrained multi-objective evolutionary algorithm, we use the multitasking-constrained multiobjective
optimization (MTCMO) framework introduced in~\cite{Qiao_etal_23}. Similar to CSO, MTCMO uses a two-population approach, and it searches for both feasible solutions and promising infeasible solutions 
by evolving a full set of candidate solutions simultaneously. This is achieved
by applying evolutionary multitasking~\cite{Gupta_etal_15} and dynamically creating auxiliary tasks to help solve the main task, i.e., solving the constrained problem itself.

Once a Pareto front has been produced, a decision-maker can investigate the alternatives and eventually select a single solution for implementation. There exists a wide range of methods for decision-making, and there even exists a range of different ways to categorize them into larger groups. They can for instance be classified as members of one of the three groups: synthesizing multi-attribute criteria methods based on utility functions, outranking methods, and interactive methods. An extensive review of these methods together with an overview of the general challenges and possibilities associated with multi-criteria decision-making is provided in~\cite{Guitouni_Martel_98}. These challenges are beyond the scope of the current paper, and we only briefly discuss the selection of Pareto solutions in the numerical results.

%\subsection{Candidate solution selection}
%{\color{blue}[Add a few references on Pareto front selection. (Maximizing social welfare, etc.)]}
%{\color{blue}[Not sure this will be discussed, or at least not have its own subsection.]}

\section{Numerical Results}
%{\color{blue}Brief intro (one-two sentences) stating what we aim to do and the needs. Repeated Multi-objective optimization in Co2 context, realistic physical model of large-scale injection site.}
In order to demonstrate the proposed multi-agent model with numerically optimized multi-objective functions available for a decision maker to select a desired strategy, we need to be able to perform repeated numerical simulation of a suitable pressure-dominated prospective CO$_2$ storage site.
As a relevant numerical test case, we consider the Bjarmeland formation located in the Barents Sea, with a setup very similar to the one described in~\cite{Allen_etal_17} and using the same open-source solver from the MATLAB Reservoir Simulation Toolbox (MRST)~\cite{Lie_2019} to facilitate comparison. The MRST CO$_2$ toolbox contains tools for representing wells, grids for realistic formations in the North Sea, various solvers for migration, pressure buildup and geomechanics, as well as a number of analysis and plotting tools~\cite{mrst_co2lab}. The Bjarmeland formation is a suitable test case in the current paper since injection is dominated by pressure buildup and maximum pressure is attained during the injection stage. Hence, the performance of different agents (wells) are expected to affect each other, and it is sufficient to consider the injection phase only. The long-term fate of the CO$_2$ plume, while essential in its own right, has no impact on the objective functions in this work and will henceforth not be further investigated. For that, and details of the setup, we refer to~\cite{Allen_etal_17} and the open-source test cases in MRST, c.f.\ \url{https://www.sintef.no/projectweb/mrst/modules/co2lab/}. A fully implicit vertically integrated black-oil type model with two phases (CO$_2$ and brine) is used to compute the migration of CO$_2$ in the reservoir~\cite{Nilsen_etal_16}. Appropriate simplifications are employed, including coarse grids and uniform rock properties, as described in~\cite{Allen_etal_17}.

We consider three injection wells located at the peaks of some of the largest structural traps of the Bjarmeland formation.  The same locations including a fourth well were investigated in~\cite{Allen_etal_17} but one well allowed only limited injection and was hence omitted in the current work. We seek the  time-varying injection rates that maximize the total amount of CO$_2$ for every coalition in a particular coalition structure. The injection period is 15 years for all wells, and is subject to change every three years. With five injection intervals per well, there is a total of 15 design variables for each optimization problem. %{\color{blue}[Comment: CSO is more ideal for large problems (many design variables) so having, say, 50 years of injection with 10 injection intervals, and four wells would probably make a better motivation for using CSO instead of Sintef's classical gradient-based methods.]} 
For simplicity, we assume that all wells start injection simultaneously, that they can change their injection rates at the same predefined times, and are subject to the same constraints imposed due to supply of CO$_2$ and minimum injection rates. We assume  minimum and maximum injection rates of respectively 0.24 and 7 Mton/year, where the former is supposed to represent a constraint based on economic feasibility and the latter a supply constraint.  The maximum rates are however intentionally chosen to be less restrictive than the physical pressure buildup constraint, defined as 90 \% of the overburden pressure and evaluated pointwise in all grid cells. The problem would otherwise simplify to the degree that simulation of any physical model would be obsolete. It is a straightforward extension to allow minimum and maximum rates as well as start and end injection times to vary between the wells, and we do not expect that to qualitatively change the results. 

%{\color{blue}[Could add explicit description of the five CS for the four agent case.]}
For constrained multi-objective and single-objective optimization, we use evolutionary methods from the MATLAB Platform for Evolutionary Multi-objective Optimization (PlatEMO)~\cite{Tian_etal_17}. More specifically, we present results for comparison using both the WSM with competitive swarm optimization (CSO) for each of the resulting SOO problems (one for each weight), and MOO using MTCMO.
Both optimization methods take as inputs the number of individuals of the population, denoted N, and the total number of model evaluations (ME). 

\begin{figure}[H]
    \centering  
\subfigure[Coalition structure: \{W1, W2\}, \{W3\}]
{\includegraphics[width=0.22\textwidth]{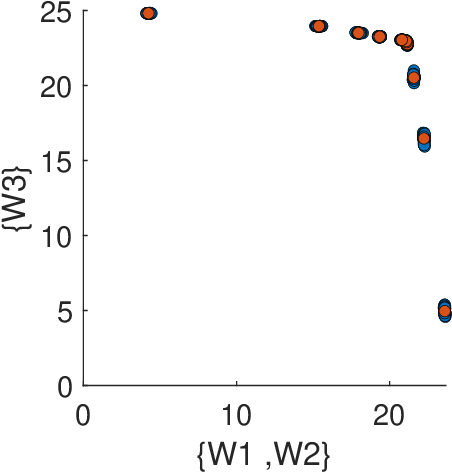}}   
~
\subfigure[Coalition structure: \{W1, W3\}, \{W2\}]
{\includegraphics[width=0.22\textwidth]{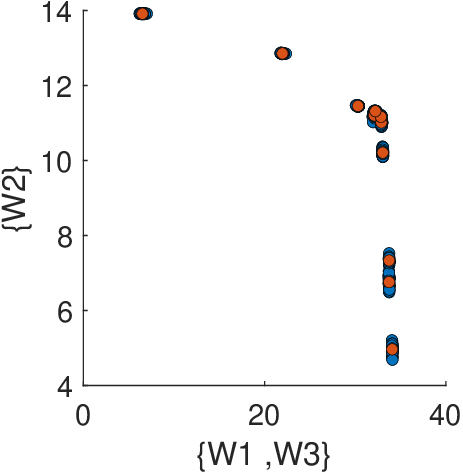}}
~
\subfigure[Coalition structure: \{W1\}, \{W2, W3\}]
{\includegraphics[width=0.22\textwidth]{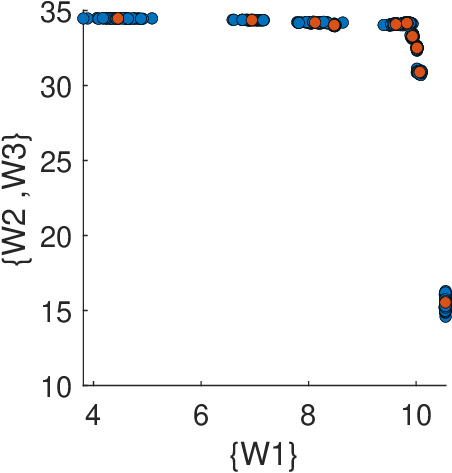}} 
~
\subfigure[Coalition structure: \{W1\}, \{W2\}, \{W3\} ]
{\includegraphics[width=0.22\textwidth]{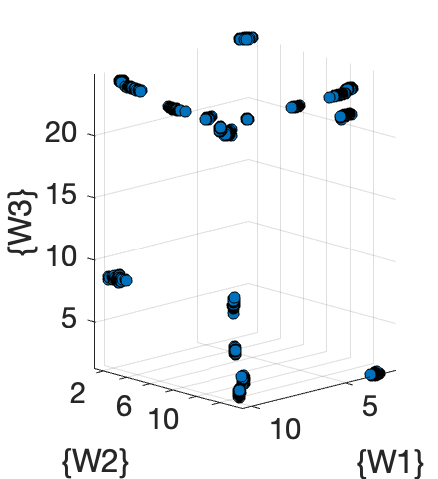}}
    \caption{Pareto fronts for the coalition structures with at least two coalitions, WSM and SOO.} %
    \label{fig:Pareto_3wells_wsm}
\end{figure}

The WSM-CSO results are shown for the Pareto fronts in Figure~\ref{fig:Pareto_3wells_wsm}, where a population of size $N=50$, and 100 iterations for a total of 5000 MRST model evaluations are used for every weight. The weights summing to 1 are equidistantly distributed with increment 0.1 for the 2D Pareto fronts and increment 0.2 for the 3D Pareto front, where it should be noted that we have also included the weight zero. 
The explicit setting of weights in the WSM method results in diverse parts of the Pareto front being captured, although the cost grows with the number of weights. The means of the 50 population members for each of the weights are shown in red (hidden in the cloud of data for the 3D Pareto front). For all cases, the Pareto fronts exhibit a characteristic kink where all coalitions receive close to their maximum attainable total injections. Note that the kink is not a numerical artifact but reflecting the characteristics of the underlying physical problem as modeled in the vertical-equilibrium solver. The Pareto front approximation method cannot capture non-convex parts, but those parts of the front that are captured can be assumed to be representative of the true front based on the observed numerical convergence. For all coalition structures with at least two coalitions, it can be observed that if a coalition injects small or moderate amounts of CO$_2$ (relative to what can be observed here), this has a small or negligible impact on the other coalition(s). This is most notable in Figure~\ref{fig:Pareto_3wells_wsm} (c) and (d). The solution for the grand coalition can be obtained by SOO, and is equivalent to the WSM solution with equal weights. By construction, this solution is always on the Pareto fronts of the other coalition structures, and by definition, it is always the Pareto solution maximizing total injections. Hence, in the current setting, it is always equal to the solution at the kinks in Figure~\ref{fig:Pareto_3wells_wsm}. This observation has been verified numerically.

\begin{figure}[H]
    \centering  
\subfigure[Optimal injection rates.]
{\includegraphics[width=0.36\textwidth]{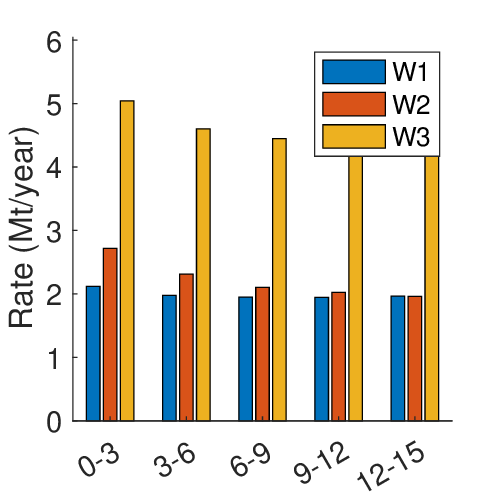}}
~
\subfigure[Plume migration, end of injection. ]
{\includegraphics[width=0.36\textwidth]{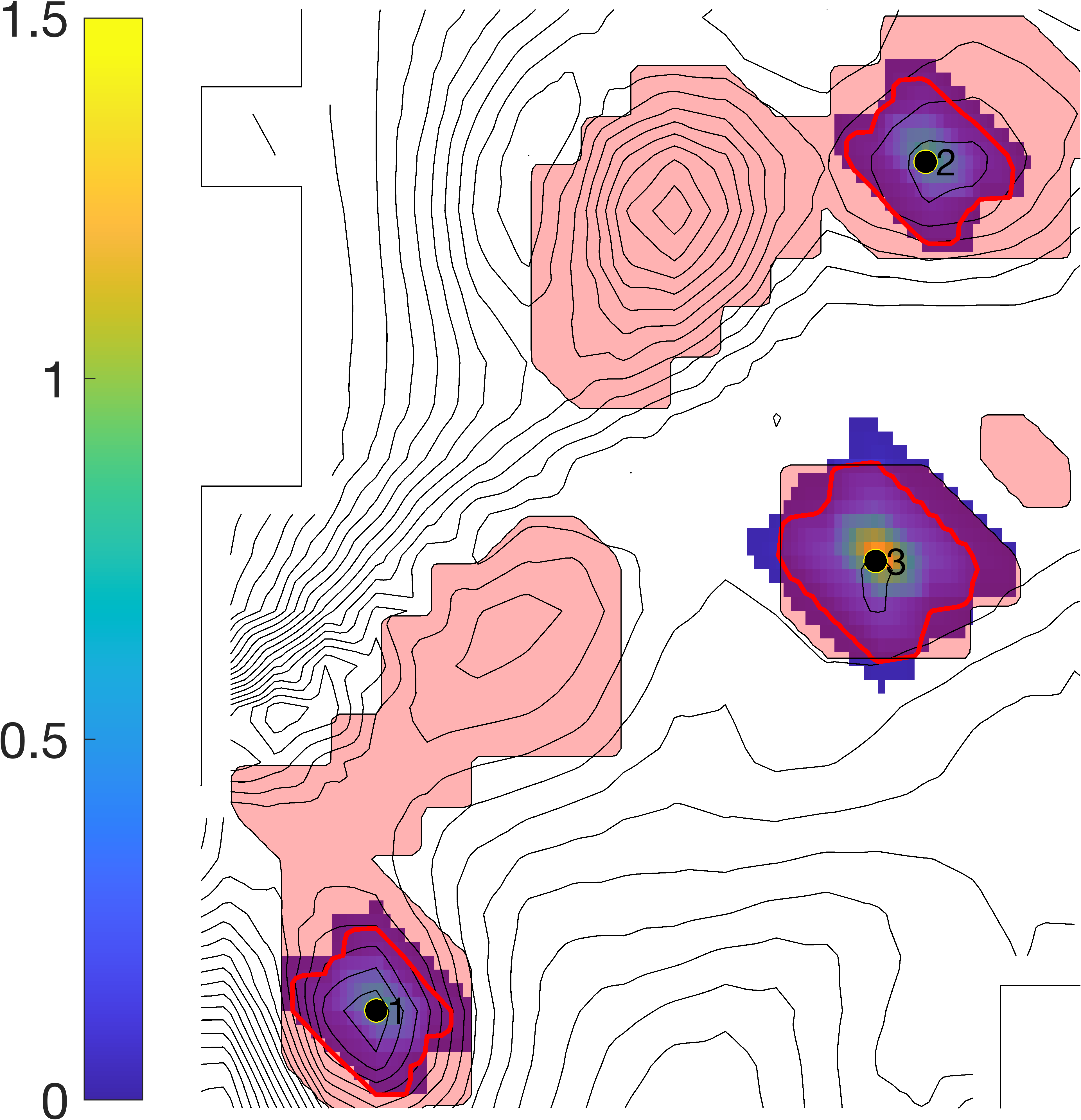}} 
~
\subfigure[Relative pressure, end of injection.]
{\includegraphics[width=0.22\textwidth]{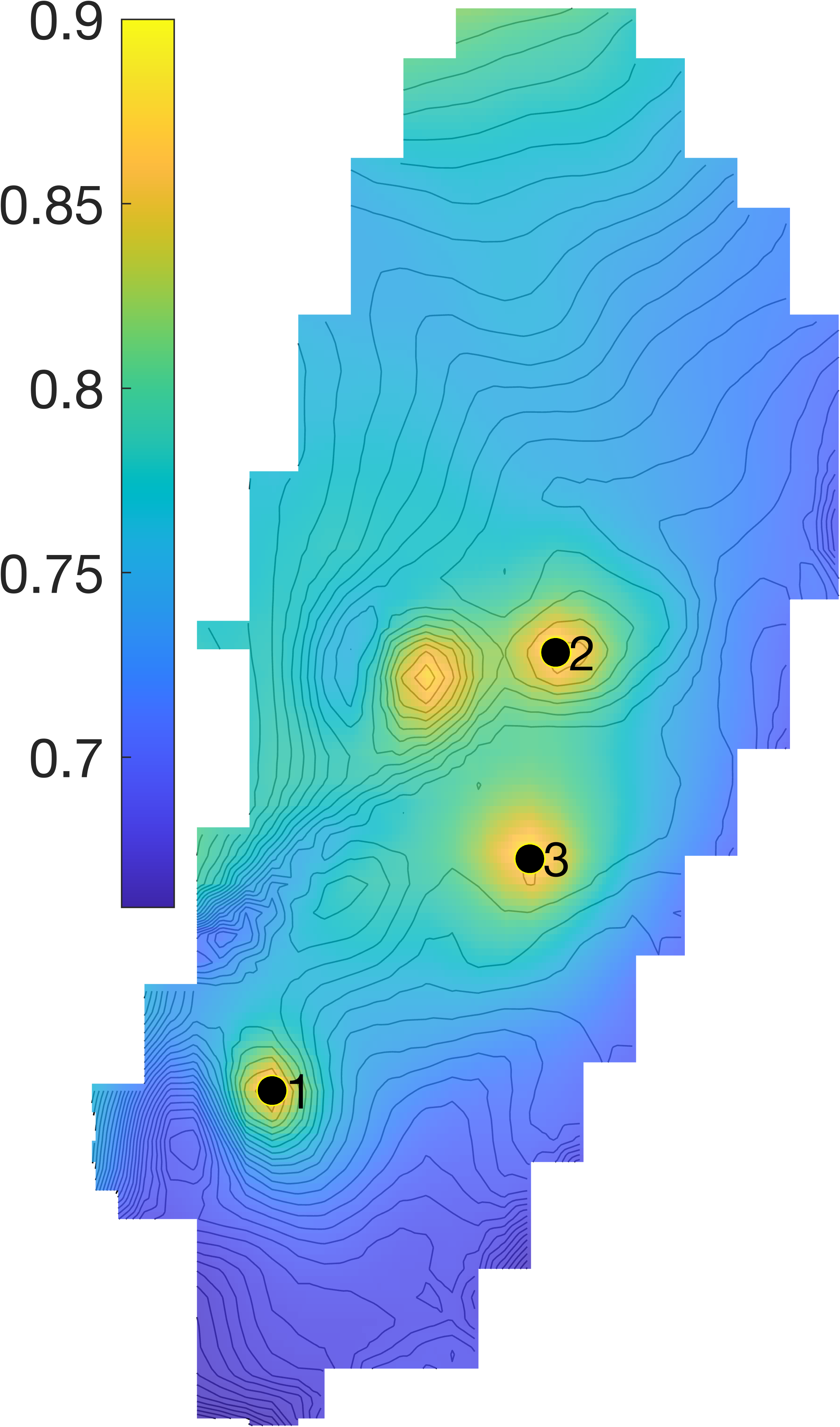}}
    \caption{Annual injection rates for all time intervals (in years), plume extension, and relative pressure buildup for the Pareto solution maximizing the total, obtained using WSM and SOO. A total of 44 Mt CO$_2$ is injected over 15 years. The color bar for the plume migration represents the mass of CO$_2$ in tons per lateral square meter, and the red contour represents 0.3 tons per m$^2$. } %
    \label{fig:All_3wells_wsm_tot_max}
\end{figure}

%{\color{blue}[Comment on the shapes of the Pareto fromnts.]}
Given the Pareto fronts of possible scenarios reflecting the goals of the coalitions, a decision-maker (e.g., some kind of legislative power regulating the injection activities without being able to fully control them) can now by means of suitable incentives enforce a selection procedure. Here we will consider two selections of Pareto solutions.
Selecting the Pareto solution that maximizes total injections (i.e., the solution at the kink) leads to almost identical optimal injection rates for all coalition structures, as expected with the linear multi-objectives. The optimal rates, plume migration, and relative pressure buildup at the end of the injection period are shown in Figure~\ref{fig:All_3wells_wsm_tot_max} for the singleton coalition structure \{\{W1\}, \{W2\}, \{W3\}\}. Since the results are very similar for all other coalition structures, they are not included here. The total amount of CO$_2$ injected over 15 years is 44 Mt. The optimized injection rates vary between wells as expected, but the variation between injection intervals is limited for all wells. For this test case, one could therefore consider having fewer injection intervals to reduce the complexity of the optimization problems. Perhaps more importantly, from a decision maker's perspective, it does not matter much what collaboration coalitions are formed since the results are so similar between coalition structures. Note that this conclusion only pertains to the situation where the Pareto solution maximizing total injections is selected. 

%{\color{blue}[Comment: This is actually something that one could prove theoretically, even for this case where we have a possibly nonlinear constraint function in the form of an MRST solver. Still, the objective function is linear. Hence, one solution of the Pareto front will always be the solution to the single-objective optimization of all wells (grand coalition). In this sense, the current reservoir simulation results are not that different from the ones we got from the semi-analytical model. We should probably point this out but perhaps save the more detailed analysis and discussion of its implications for the full-length peer-review paper.]} 

\begin{figure}[H]
    \centering  
\subfigure[CS: \{W1, W2\},    \{W3\} Total amount 28 Mt.]
{\includegraphics[width=0.22\textwidth]{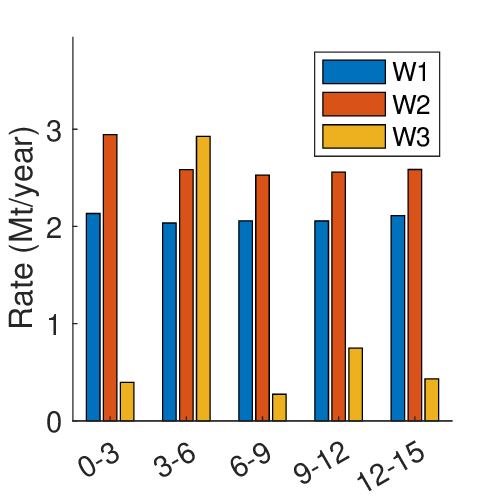}}
\subfigure[CS: \{W1, W3\}    \{W2\} Total amount 39 Mt.]
{\includegraphics[width=0.22\textwidth]{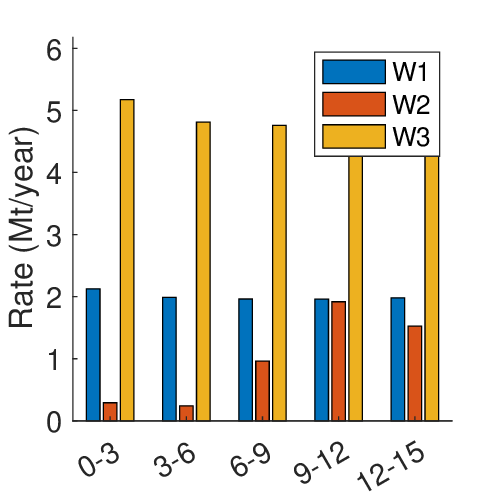}} \subfigure[CS: \{W1\},    \{W2, W3\} Total amount 26 Mt.]
{\includegraphics[width=0.22\textwidth]{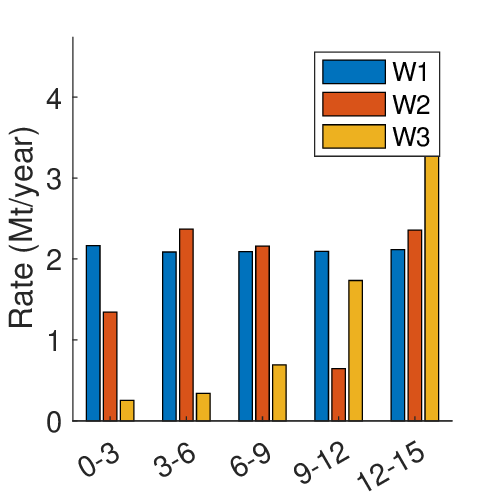}}
\subfigure[CS:~\{W1\},\{W2\},\{W3\} Total amount 24 Mt.]
{\includegraphics[width=0.22\textwidth]{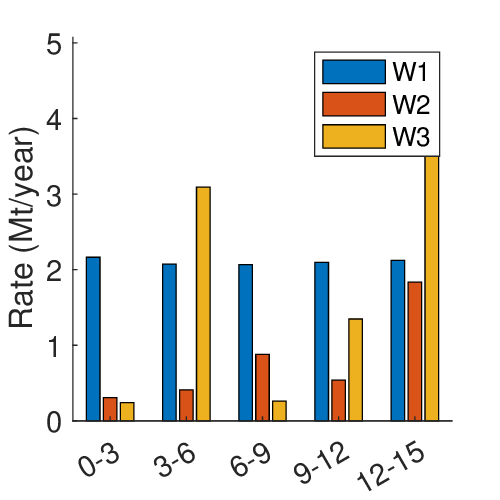}}
    \caption{Injection rates for the Pareto solution maximizing the injections of W1 for different coalition structures (CS), obtained using WSM and SOO.} %
    \label{fig:Rates_3wells_wsm_maxW1}
\end{figure}

\begin{figure}[H]
    \centering  
\subfigure[CS: \{W1, W2\}, \{W3\}]
{\includegraphics[width=0.22\textwidth]{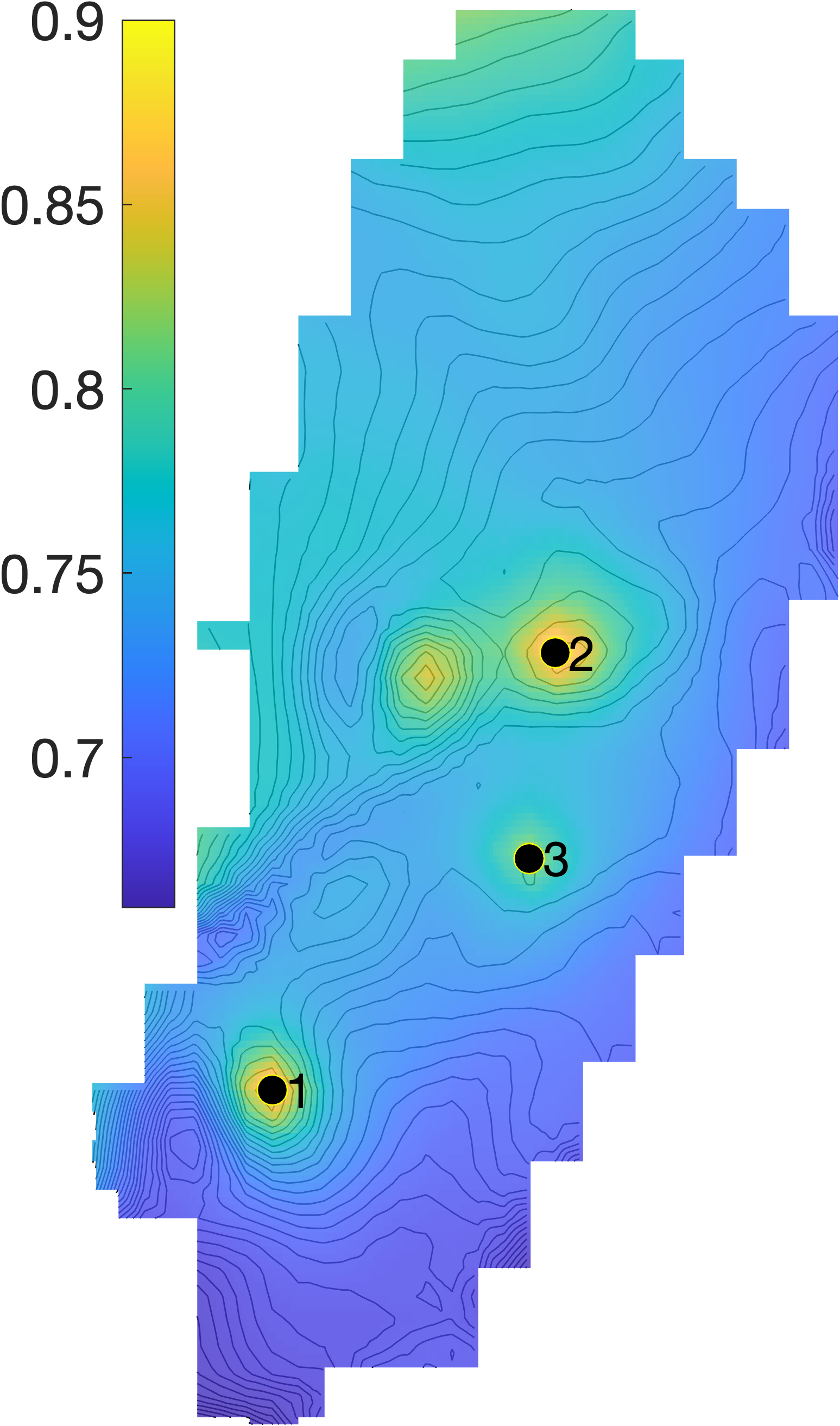}}
\subfigure[CS: \{W1, W3\}, \{W2\} ]
{\includegraphics[width=0.22\textwidth]{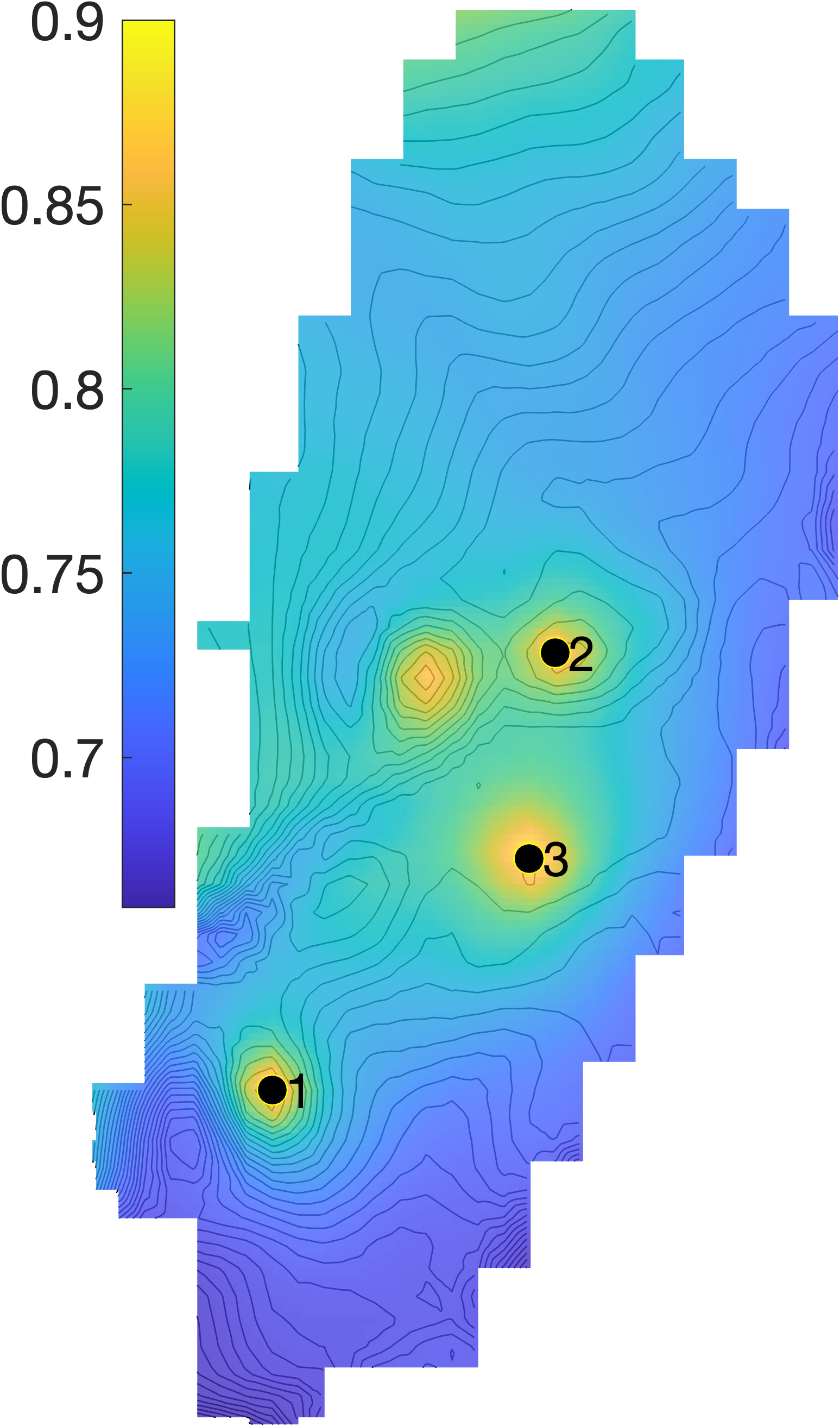}} \subfigure[CS: \{W1\},    \{W2, W3\}]
{\includegraphics[width=0.22\textwidth]{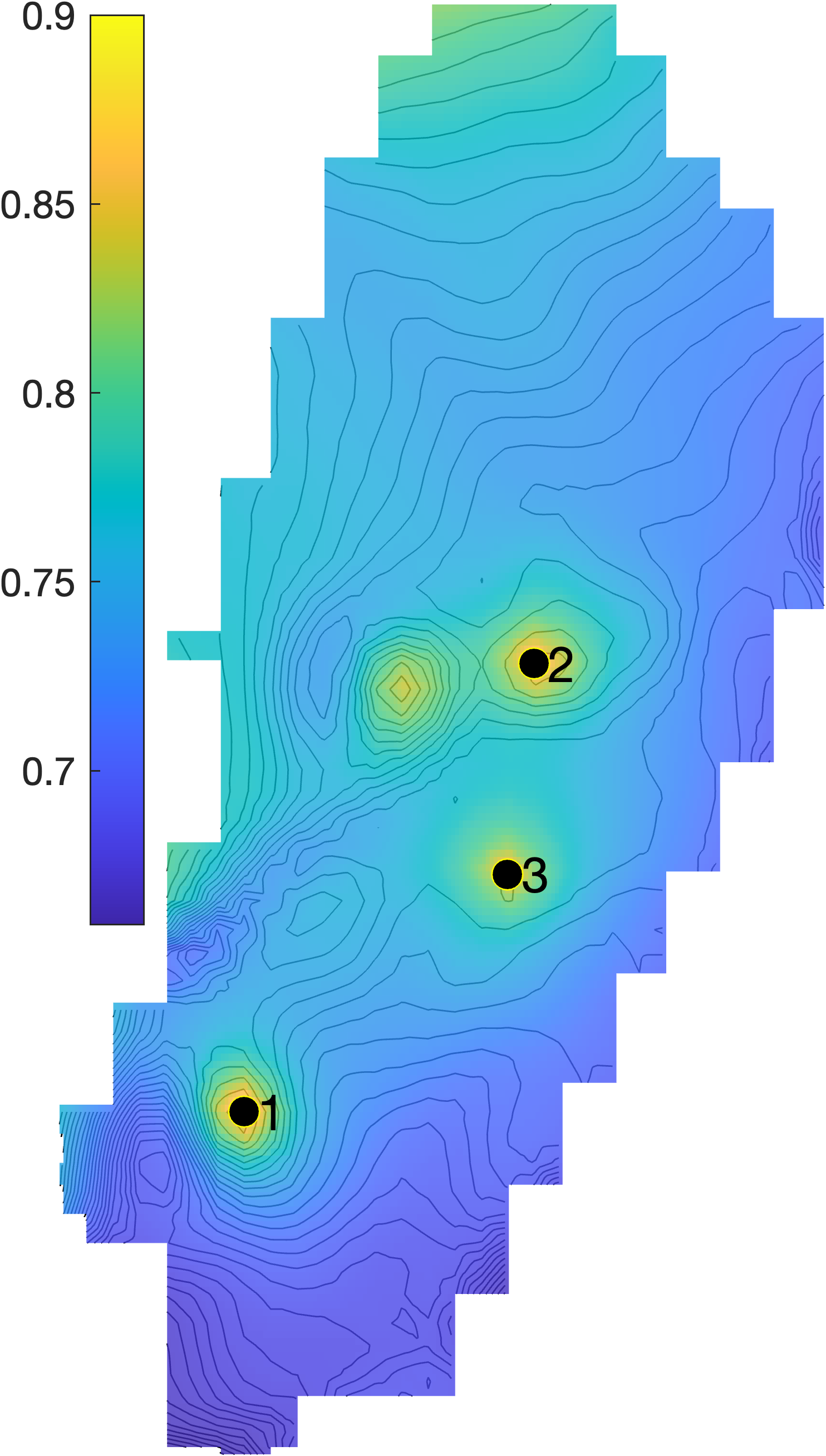}}
\subfigure[CS: \{W1\},\{W2\},\{W3\}]
{\includegraphics[width=0.22\textwidth]{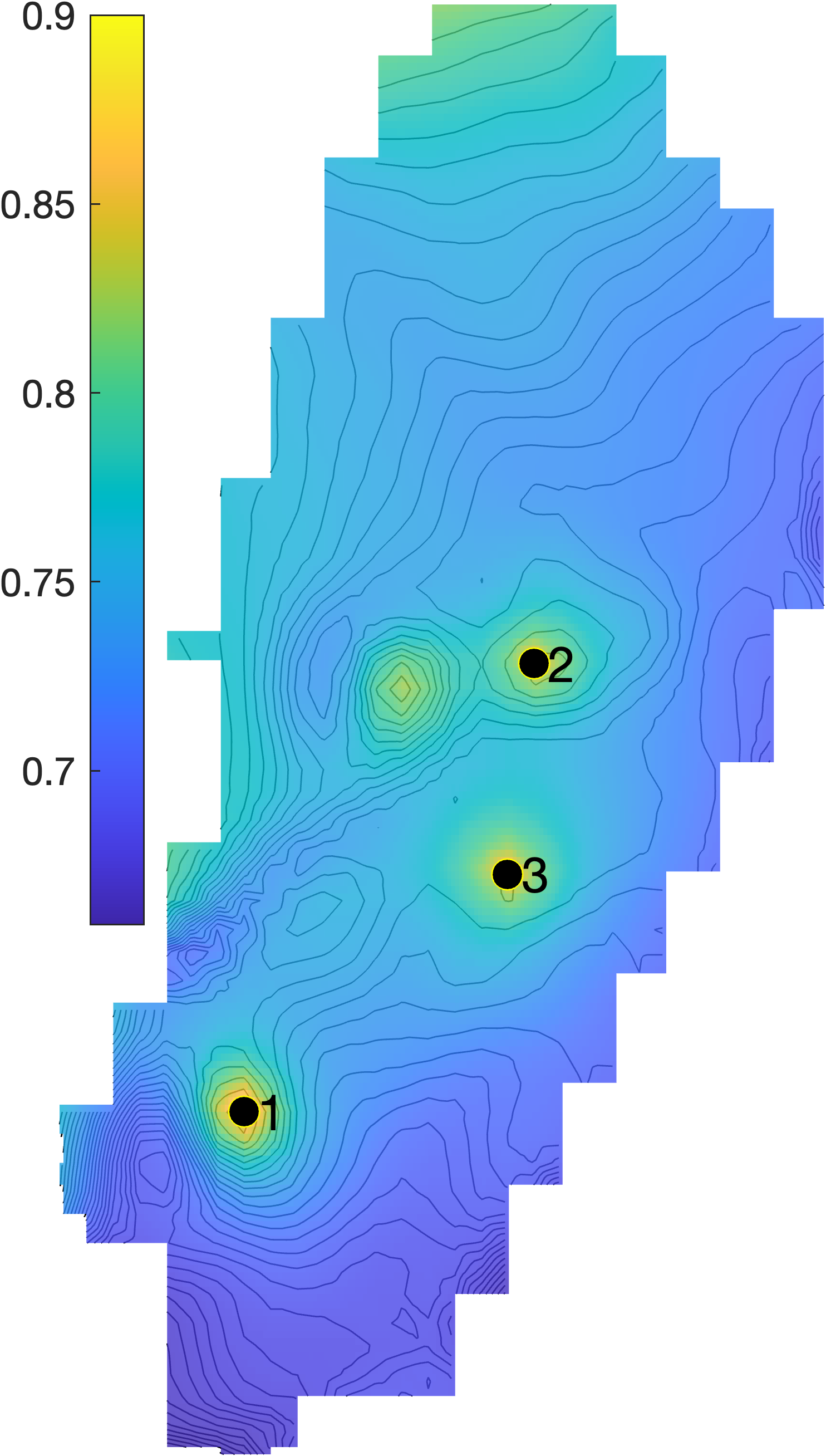}}
    \caption{Relative pressure with respect to overburden pressure at the end of the injection period. Pareto solution maximizing the injections of W1 for different coalition structures, obtained using WSM and SOO.} %
    \label{fig:Pressure_3wells_wsm_maxW1}
\end{figure}

\begin{figure}[H]
    \centering  
\subfigure[CS: \{W1, W2\}, \{W3\}]
{\includegraphics[width=0.225\textwidth]{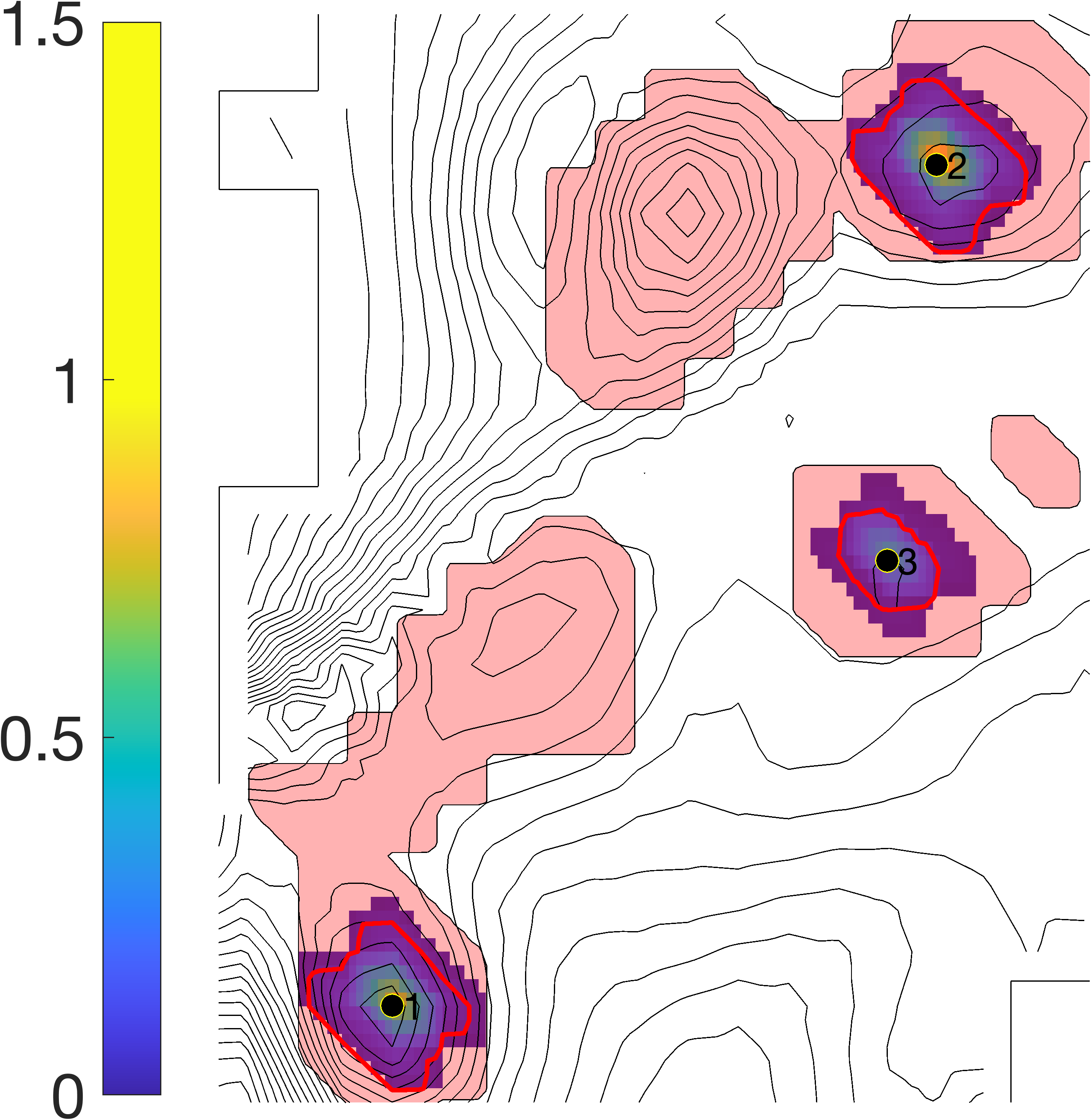}}
\subfigure[CS: \{W1, W3\}, \{W2\} ]
{\includegraphics[width=0.225\textwidth]{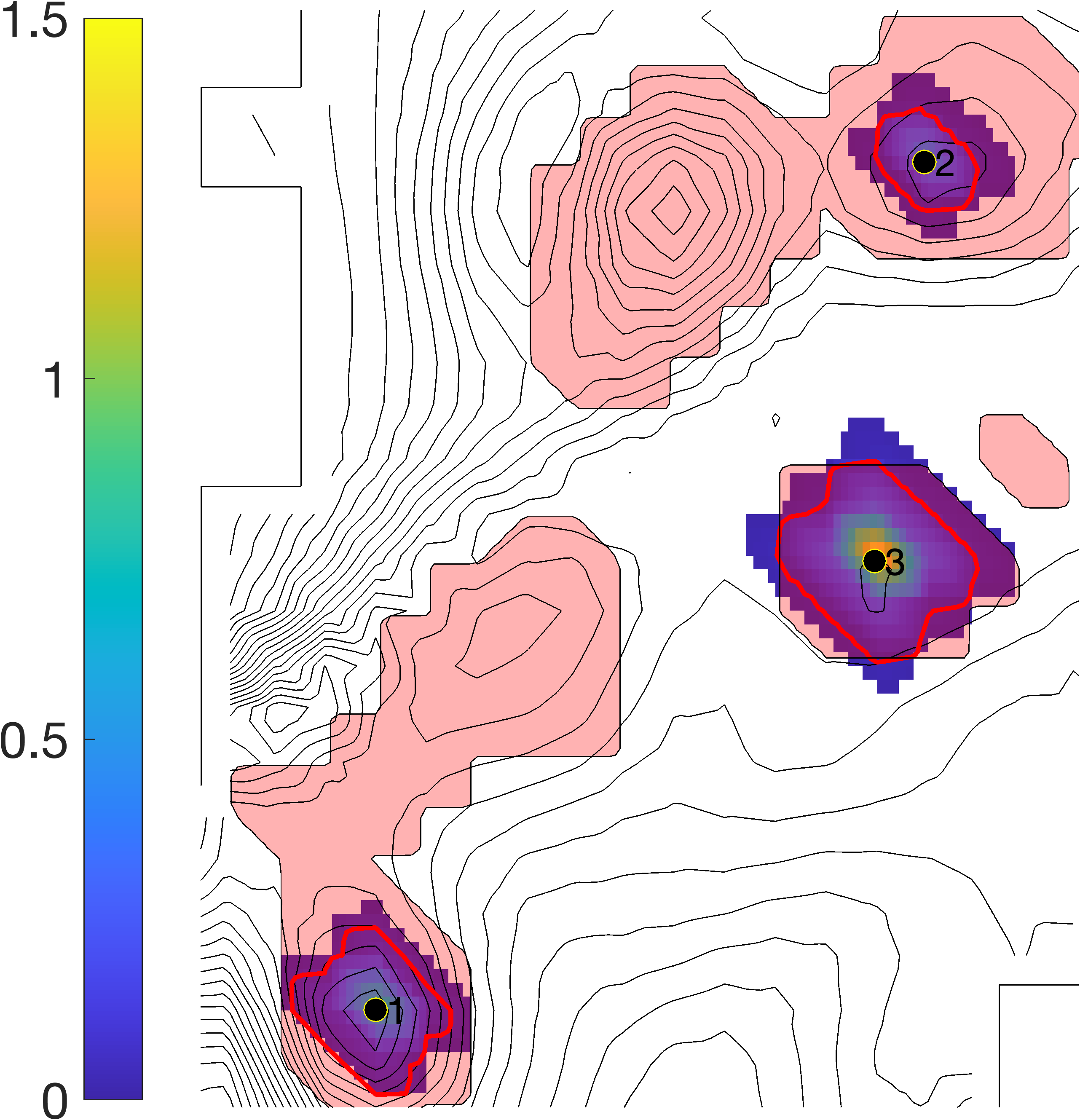}} \subfigure[CS: \{W1\},    \{W2, W3\}]
{\includegraphics[width=0.225\textwidth]{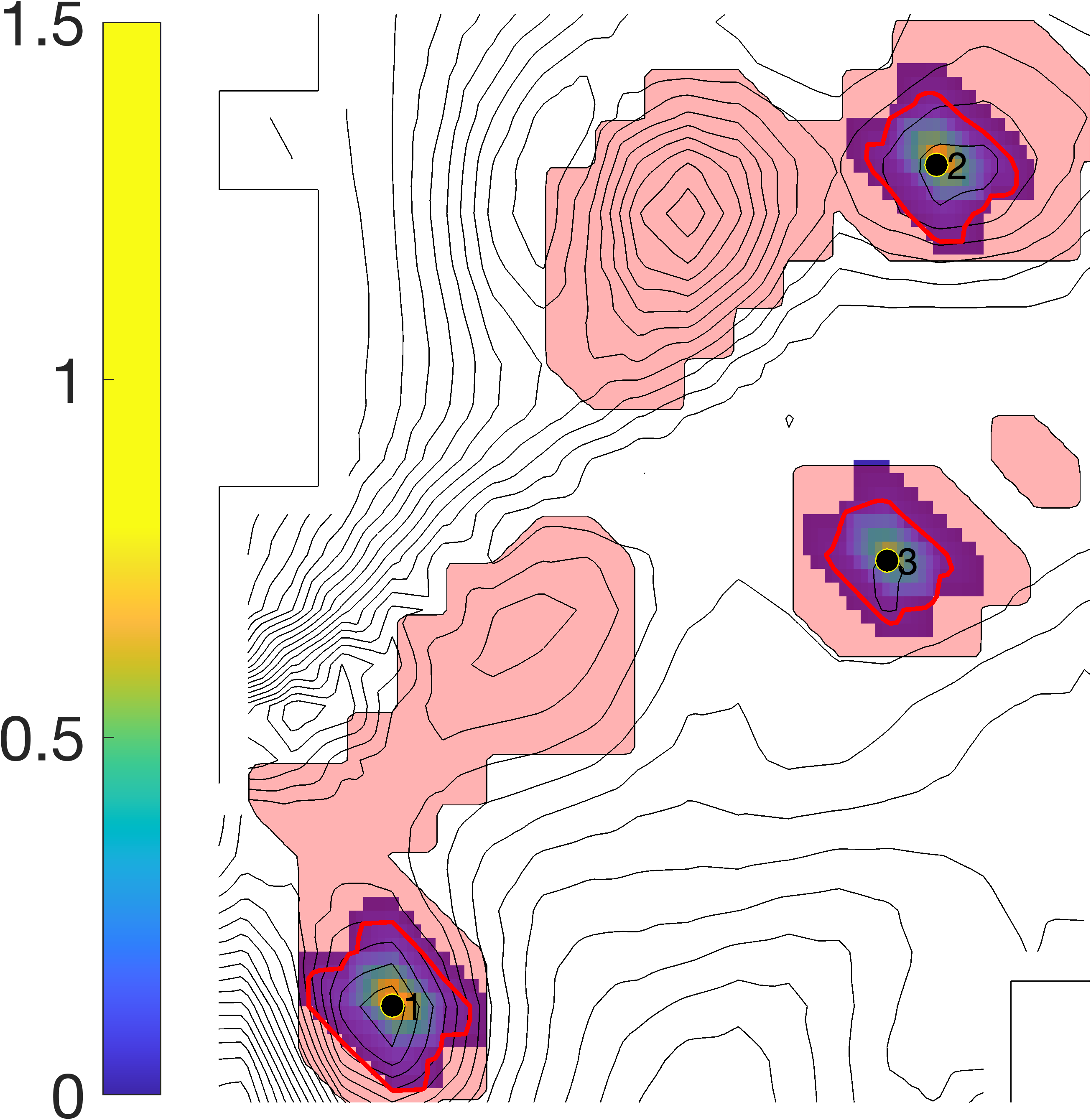}}
\subfigure[CS: \{W1\},\{W2\},\{W3\}]
{\includegraphics[width=0.23\textwidth]{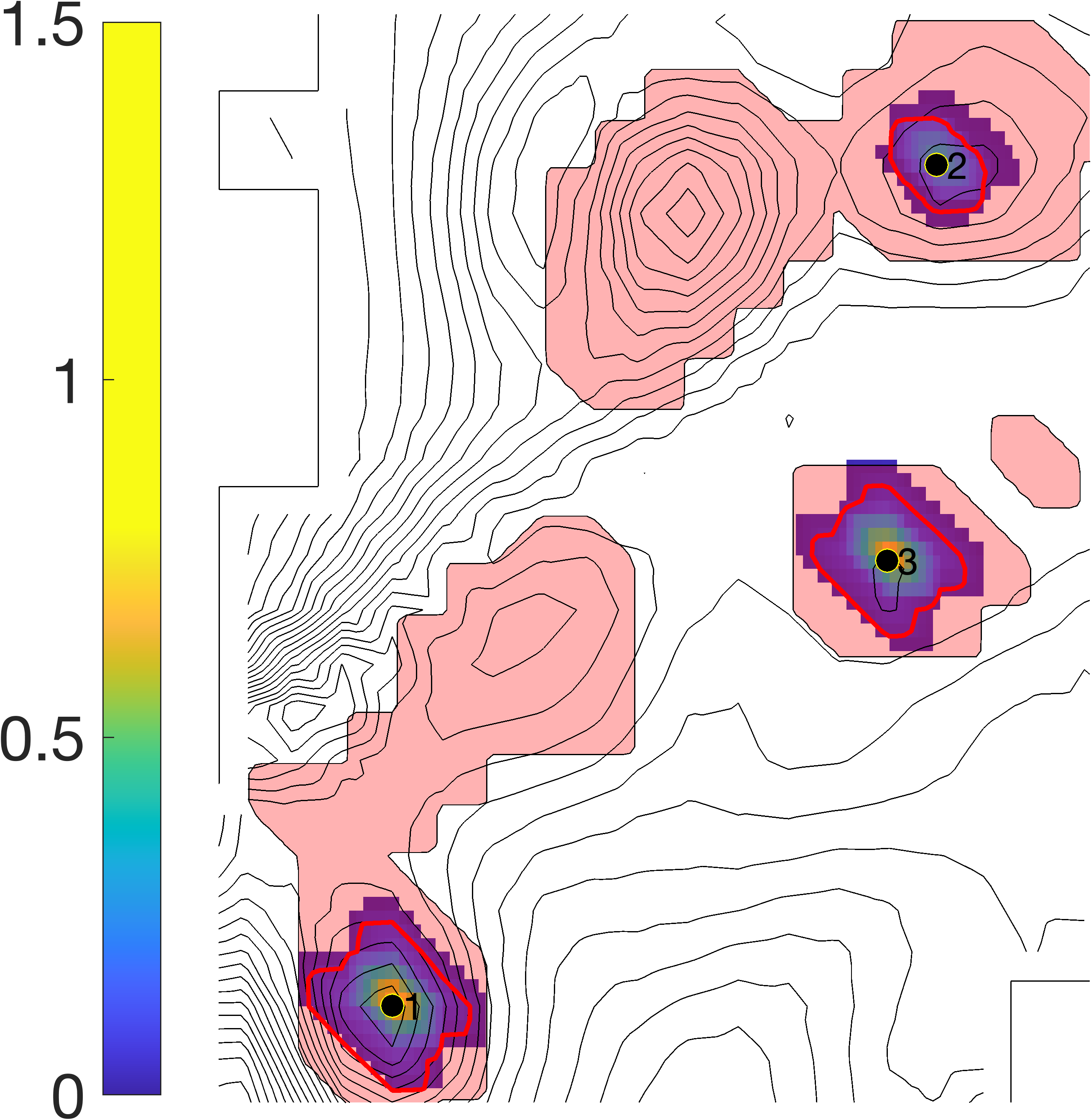}}
    \caption{CO$_2$ plume at the end of the injection period. Pareto solution maximizing the injections of W1, obtained using WSM and SOO. The color bar for the plume migration represents mass of CO$_2$ in tons per lateral square meter, and the red contour represents 0.3 tons per m$^2$.} %
    \label{fig:Plume_3wells_wsm_maxW1}
\end{figure}

As an alternative Pareto front selection criterion that leads to more variability in the solutions, we next assume that one wants to favor solutions where the injections of W1 are as large as possible (whether W1 is in a coalition with others or by itself). The resulting injection rates for the four coalition structures excluding the grand coalition are shown in Figure~\ref{fig:Rates_3wells_wsm_maxW1}. The rates of the favored well W1 do not vary much, but the well being in the same coalition as W1 indeed benefits from that. %{\color{blue}[Not sure everything here is quite trustworthy, or at least I don't know why the rates are varying so much between injection intervals in some cases. Either the results are noisy/not converged and/or there is some similar effect making them overly sensitive.]} 
The corresponding relative pressure buildup in the Bjarmeland formation is shown in Figure~\ref{fig:Pressure_3wells_wsm_maxW1}, reflecting the behavior seen in the injection rates. The CO$_2$ plumes are depicted in Figure~\ref{fig:Plume_3wells_wsm_maxW1}, indicating that most of the injected CO$_2$ at this early stage is still confined within the structural traps.

\begin{figure}[H]
    \centering  
\subfigure[N=50, ME=5000]
{\includegraphics[width=0.22\textwidth]{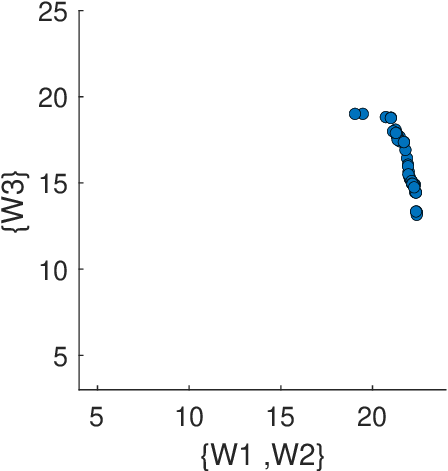}}
\subfigure[N=100, ME=10000]
{\includegraphics[width=0.22\textwidth]{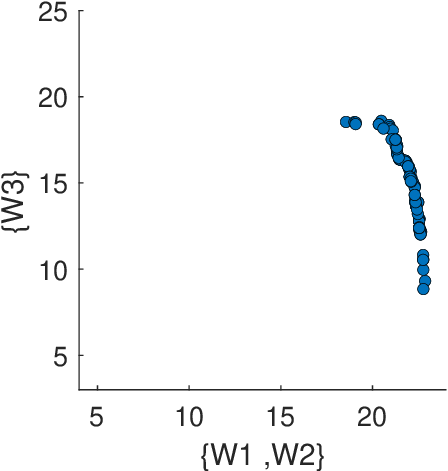}} \subfigure[N=200, ME=10000]
{\includegraphics[width=0.22\textwidth]{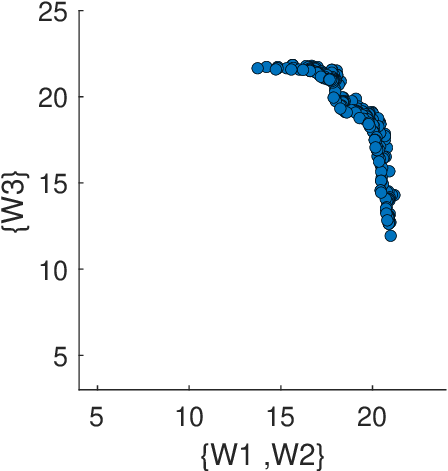}}
\subfigure[N=500, ME=25000]
{\includegraphics[width=0.22\textwidth]{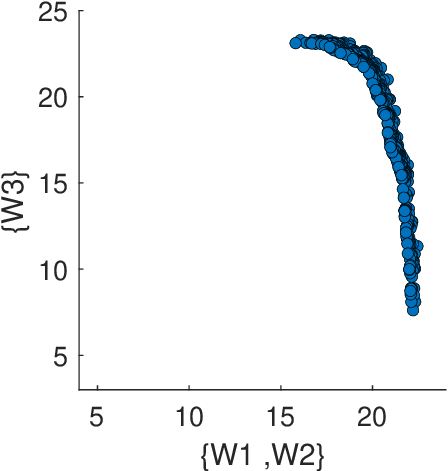}}
    \caption{Pareto fronts for the total amount of CO$_2$ (all values in Mt) with coalition structure \{\{W1,W2\},\{W3\}\}, obtained with MTCMO with varying population size (N) and number of physical model evaluations (ME).} %
    \label{fig:Pareto_3wells_MTCMO_W12_W3}
\end{figure}

\begin{figure}[H]
    \centering  
\subfigure[N=50, ME=5000]
{\includegraphics[width=0.22\textwidth]{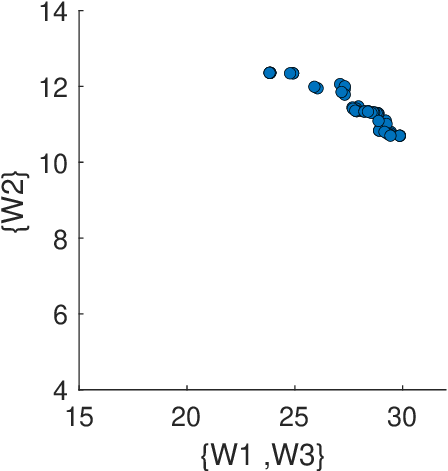}}
\subfigure[N=100, ME=10000]
{\includegraphics[width=0.22\textwidth]{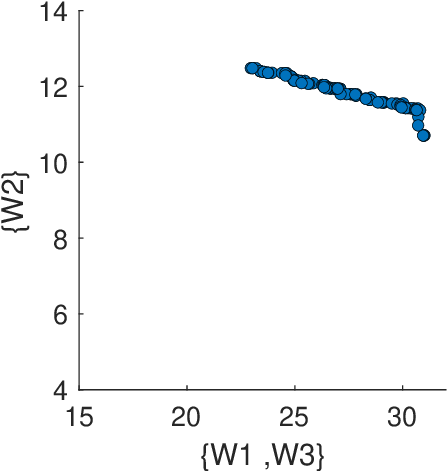}} \subfigure[N=200, ME=10000]
{\includegraphics[width=0.22\textwidth]{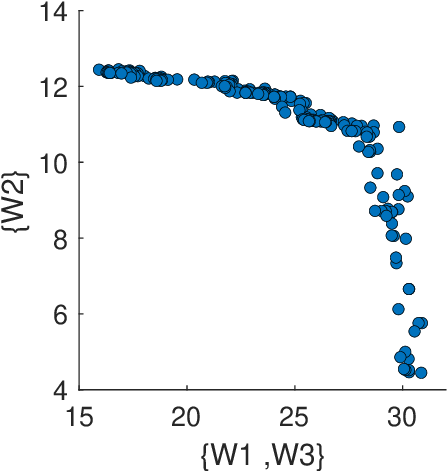}}
\subfigure[N=500, ME=25000]
{\includegraphics[width=0.22\textwidth]{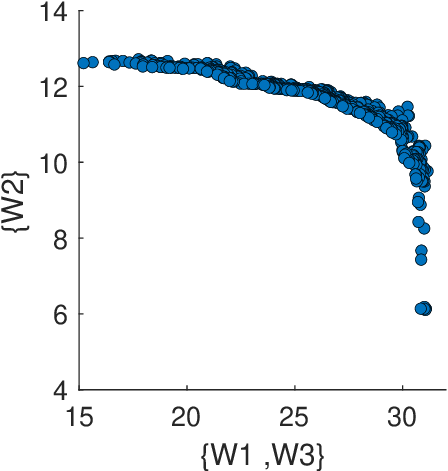}}
    \caption{Pareto fronts for the total amount of CO$_2$ (Mt) with coalition structure \{\{W1,W3\},\{W2\}\}, obtained with MTCMO with varying population size (N) and number of physical model evaluations (ME).} %
    \label{fig:Pareto_3wells_MTCMO_W13_W2}
\end{figure}

\begin{figure}[H]
    \centering  
\subfigure[N=50, ME=5000]
{\includegraphics[width=0.22\textwidth]{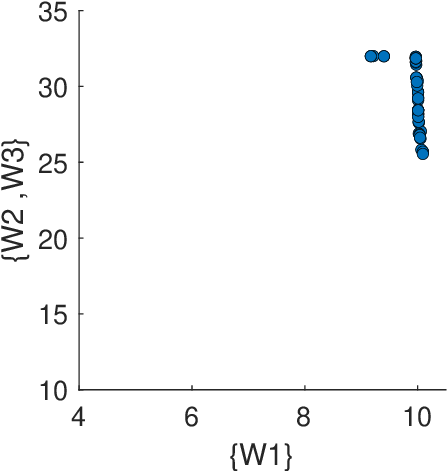}}
\subfigure[N=100, ME=10000]
{\includegraphics[width=0.22\textwidth]{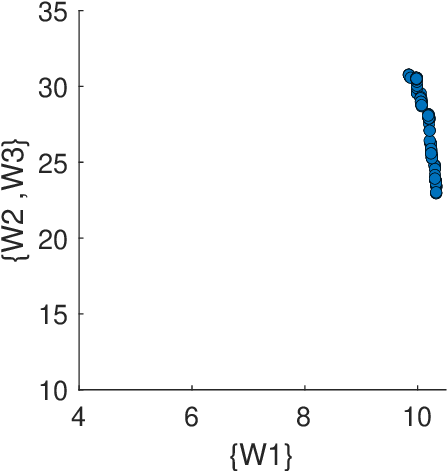}} \subfigure[N=200, ME=10000]
{\includegraphics[width=0.22\textwidth]{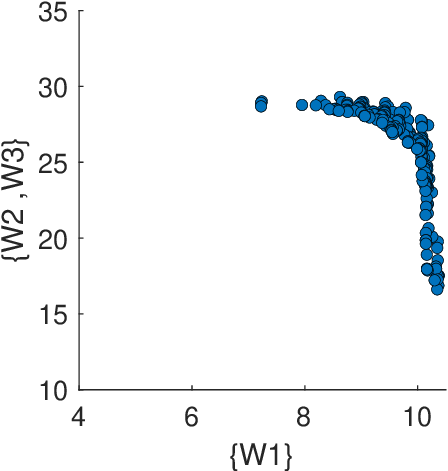}}
\subfigure[N=500, ME=25000]
{\includegraphics[width=0.22\textwidth]{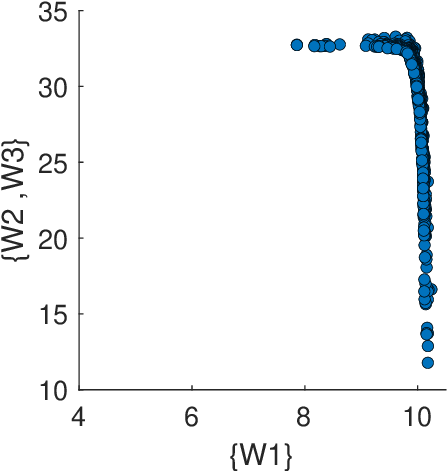}}
    \caption{Pareto fronts for the total amount of CO$_2$ (in Mt) with coalition structure \{\{W1\},\{W2,W3\}\}, obtained with MTCMO with varying population size (N) and number of physical model evaluations (ME).} %
    \label{fig:Pareto_3wells_MTCMO_W1_W23}
\end{figure}

\begin{figure}[H]
    \centering  
\subfigure[N=50, ME=5000]
{\includegraphics[width=0.22\textwidth]{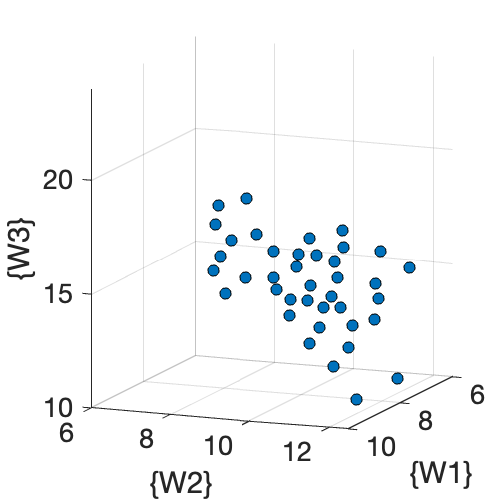}}
\subfigure[N=100, ME=10000]
{\includegraphics[width=0.22\textwidth]{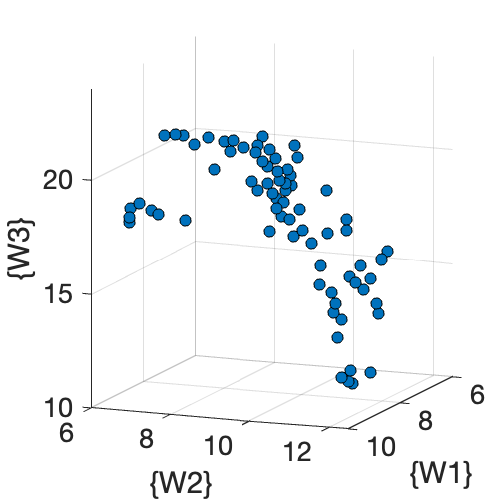}} \subfigure[N=200, ME=10000]
{\includegraphics[width=0.22\textwidth]{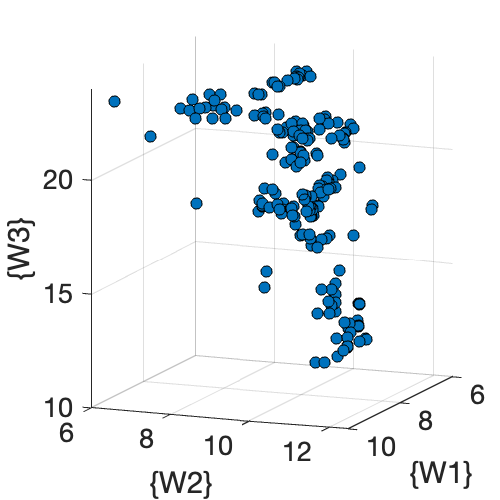}}
\subfigure[N=500, ME=25000]
{\includegraphics[width=0.22\textwidth]{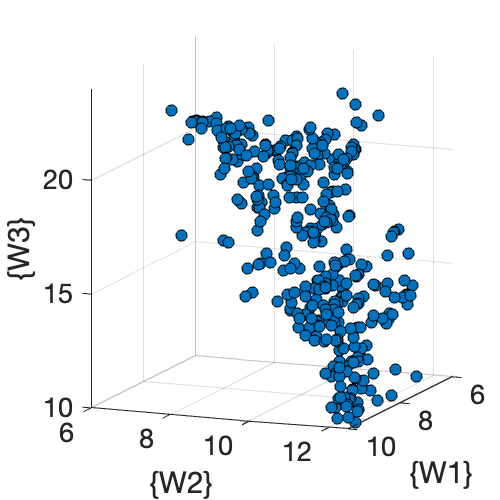}}
    \caption{Pareto fronts for the total amount of CO$_2$ (in Mt) with coalition structure \{\{W1\},\{W2\},\{W3\}\}, obtained with MTCMO with varying population size (N) and number of physical model evaluations (ME). } %
    \label{fig:Pareto_3wells_MTCMO_W1_W2_W3}
\end{figure}

\begin{figure}[H]
    \centering  
\subfigure[CS: \{W1,W2\},\{W3\}]
{\includegraphics[width=0.225\textwidth]{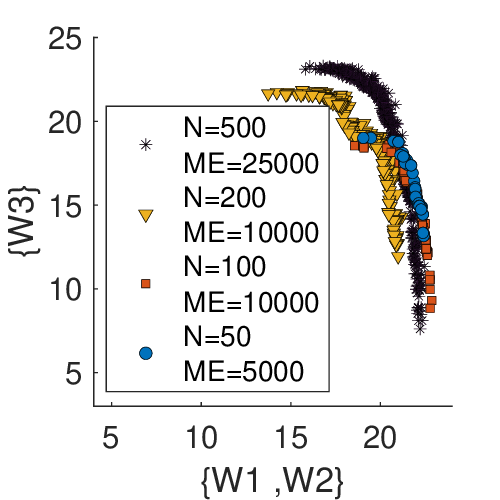}}
\subfigure[CS: \{W1,W3\},\{W2\}]
{\includegraphics[width=0.225\textwidth]{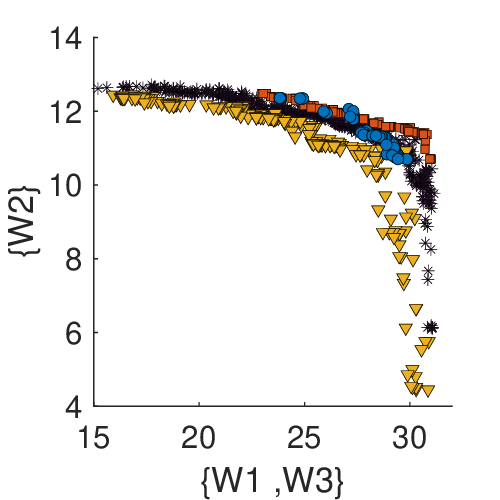}} 
\subfigure[CS: \{W1\},\{W2,W3\}]
{\includegraphics[width=0.225\textwidth]{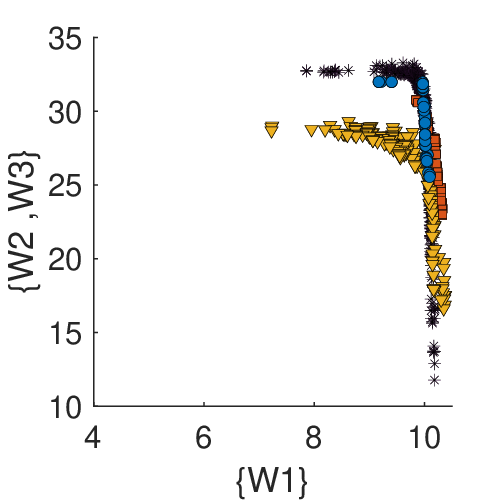}}
\subfigure[CS:~\{W1\},\{W2\},\{W3\}]
{\includegraphics[width=0.229\textwidth]{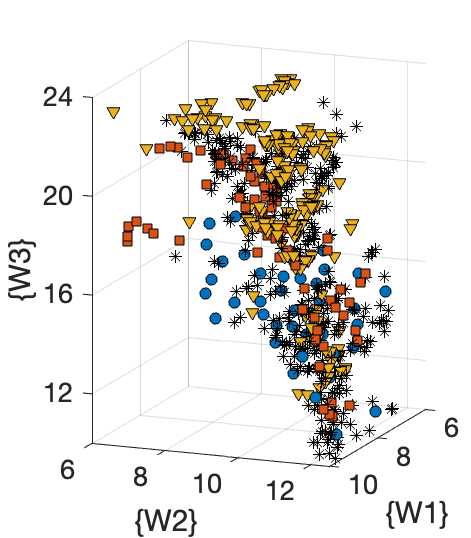}}
    \caption{Summary of all Pareto fronts for total CO$_2$ injections (in Mt), obtained with MTCMO with varying population size (N) and number of physical model evaluations (ME). The legends in (a) pertain to all subfigures.} %
    \label{fig:Pareto_3wells_MTCMO_summary}
\end{figure}

Next, we compare the Pareto fronts obtained using the WSM with truly multi-objective optimization. The same setup as before is used and solved with variable population sizes and maximum number of function evaluations using the MTCMO previously described and shown in Figures~\ref{fig:Pareto_3wells_MTCMO_W12_W3}, \ref{fig:Pareto_3wells_MTCMO_W13_W2}, \ref{fig:Pareto_3wells_MTCMO_W1_W23}, and \ref{fig:Pareto_3wells_MTCMO_W1_W2_W3}.
For smaller populations, the ends of the Pareto fronts are not captured with MTCMO. The trends are most easily observed for the 2D Pareto fronts in Figures~\ref{fig:Pareto_3wells_MTCMO_W12_W3}, \ref{fig:Pareto_3wells_MTCMO_W13_W2}, and~\ref{fig:Pareto_3wells_MTCMO_W1_W23}.
The 3D Pareto fronts in Figure~\ref{fig:Pareto_3wells_MTCMO_W1_W2_W3} are included for completeness and show similar patterns, although less clearly. Finally, all the Pareto fronts from MTCMO are shown together to facilitate comparison in Figure~\ref{fig:Pareto_3wells_MTCMO_summary}.

A much larger portion of the Pareto front seems to be captured by the WSM-SOO compared to MOO using  MTCMO, despite the fact that the latter is supposed to be able to capture parts of the fronts that the former cannot capture. Other MOO methods from the PlatEMO package exhibited similar or worse behavior, and the results are not included here. With increasing number of model evaluations, and in particular with increasing population size, larger sections of the Pareto fronts are captured. The largest sets of 25000 model evaluations are smaller than the total number of model evaluations used for the WSM-SOO method (55000 model evaluations for the 2D Pareto fronts 105000 model evaluations for the 3D front), so better coverage would be expected for the MTCMO method in a comparison with identical number of model evaluations. It is worthwhile to note that the different MTCMO approximations of the same Pareto fronts differ not only in their ability to span the full front, but they also yield results that in some cases appear as distinctly different fronts, as seen in the summary Figure~\ref{fig:Pareto_3wells_MTCMO_summary}. The population members of MTCMO (and related methods) exchange information in order to accelerate convergence of the Pareto front. However, it appears that this can result in Pareto front approximations that appear to be numerically converged but in reality are relatively far from the true Pareto optimal solutions.

% %%%%%%
% \begin{figure}[H]
%     \centering  
% \subfigure[ECMO]
% {\includegraphics[width=0.22\textwidth]{figures_ECMOR/Pareto_W12_W3_MOO_EMCMO.eps}}
% \subfigure[MTCMO]
% {\includegraphics[width=0.22\textwidth]{figures_ECMOR/Pareto_W12_W3_MOO_MTCMO.eps}}
%     \caption{Pareto fronts for three injection wells, MOO. {\color{blue}[This figure produced when testing the performance of different MOO methods will be removed.]}} %
%     \label{fig:Pareto_3wells_moo}
% \end{figure}

%%%%%%%%
%%%%%%%%
%%%%%
%%%
%%%%%%

%%%
% %%%%%%
%%%%%%

%%%%%%
%%%%%%

%%%%%
% \begin{figure}[H]
%     \centering  
% \subfigure[]
% {\includegraphics[width=0.22\textwidth]{figures_ECMOR/Pressure_map_EX.png}}
% \subfigure[]
% {\includegraphics[width=0.22\textwidth]{figures_ECMOR/Plume_map_EX.png}} \subfigure[]
% {\includegraphics[width=0.44\textwidth]{figures_ECMOR/Injection_schedule_EX.eps}}

%     \caption{Example figures.} %
%     \label{fig:ext_utopia_wts_4wells}
% \end{figure}

%\textbf{Possible variations for numerical investigation}\\
%- Single-objective and WSM vs multi-objective CSO for Pareto front approximation: cost and accuracy.

\section{Discussion}
%{\color{blue}[Discussion not part of the template. Could move the text to numerical results and conclusions, respectively]}

%{\color{blue}[Can write about how we would model the situation where there are agents that operate two or more wells - that's an interesting discussion but might also be saved for later use.]}

The complexity of the multi-agent problem quickly becomes infeasible with increasing number of agents as the number of coalition structures described by the Bell number grows even faster than the exponential growth of the number of coalitions. While this needs to be considered in practical applications, there are realistic situations of CO$_2$ injection with numbers of wells much larger than what was considered in the numerical experiments, that could still be targeted with the presented method. First of all, the number of agents that control the wells is not likely to be very large. Second, not all possible coalitions need to be simulated if they are not realistic. If a certain collaboration is not of interest, the corresponding coalition structure does not need to be further investigated. Also, with a large number of wells and when there are agents that operate two or more wells, there may be relatively few coalition structures determined by the number of operators, and not the number of wells. However, the complexity of each multi-objective optimization problem corresponding to a given coalition structure grows with the number of wells, so there is certainly a cost associated with the number of wells, although it does not grow as fast as with the number of agents.

%Weighted sum method vs MO methods - pros and cons: simulation speed, accuracy etc.

The WSM-SOO is often described as a classical method in the literature, probably due to its relative age, and its simplicity to use and interpret. In the current work, it appears to capture the Pareto front despite its shortcomings in capturing non-convex Pareto fronts. The objective function is convex, so any non-convex parts - if there are any - must be imposed by the overburden pressure constraint function. With uniform and equal refinements of the weights for all coalitions, the WSM-SOO method suffers from the curse of dimensionality, i.e., the numerical cost grows exponentially with the number of coalitions. Hence, other weights need to be considered for more complex problems.

%Numerical cost of the two approaches WSM-SOO vs MOO with the number of agents (and injection control intervals, the number of which is always linear in the number of agents in the current setting). Using tensor-like grids for the weights {\color{blue}[Check whether 'tensor-like grids' is really an appropriate description]}, the WSM suffers from the curse of dimensionality in the cost of a single transformed MOO (i.e., a single coalition structure). Then, in addition, the number of coalition structures grows exponentially with the number of agents, regardless of the choice of optimization method.

%{\color{blue}[How to make a decision based on the Pareto fronts.]}
The numerical results with MOO show a clear pattern for how the Pareto front is captured, depending on the population size and the number of generations. The population size determines how well the space of potential solutions is captured, with a clear tendency to find the ends of the front only if the population size is large, while the middle part that typically exhibits a clear kink requires significantly fewer model evaluations. Hence, if one is mostly interested in the Pareto candidate solution that maximizes total efficiency, i.e., the total amount of injected CO$_2$, MOO can produce the solution to reasonable accuracy with a relatively small number of model evaluations.

There is a trade-off between the accuracy of the physical model describing the pressure buildup during injection, and the numerical cost for each model evaluation that needs to be extensively repeated.
%Is the physical model sufficiently accurate? Short simulation time a requirement, but 
Once a tentative decision has been made based on a Pareto efficient solution using simplified-physics models, a single or a small number of high-fidelity simulations can be performed for validation of the selected solution. Such a test would make sure that the expected rates obtained in the large-scale multi-objective optimization can indeed be attained, but it clearly cannot show that they are the optimal ones. This combination of many low-fidelity and a few high-fidelity simulations is expected to yield satisfactory results if the low-fidelity solutions are accurate enough to yield sufficiently close to the true optimal rates.

%{\color{blue}Need for longer injection times, more wells, improved physical model, etc.}

The total injection times considered here are relatively short, in particular for a pressure-limited injection site. To exploit the full storage potential of the formation, longer injection periods should be investigated within the multi-agent framework. The advantages of variable injection rates are likely to become more clear with longer injection times. As shown in previous work, relocation of the wells can also significantly increase the amount of CO$_2$ that can be safely stored~\cite{Allen_etal_17}.

\section{Summary and Conclusions}
Large-scale deployment of subsurface resources for CO$_2$ storage by independent agents will require injection strategies that include incorporation of physical constraints affected both by the formation properties and the activities of competing agents performing CO$_2$ or hydrogen injection, or performing hydrocarbon extraction.
We present a multi-agent model and constrained multiobjective optimization of CO$_2$ injection, where the operators of injection wells can choose to form binding agreements to collaborate to maximize their joint performance. The multi-agent model is demonstrated on the pressure-limited Bjarmeland formation in the Barents Sea. Pareto fronts depicting candidate injection schedules are obtained with both the weighted sum method that transforms multi-objective optimization problems to a set of single-objective optimization problems, and with an evolutionary two-population method for truly multi-objective optimization. For simplicity, an agent is supposed to represent a single CO$_2$ injection well. The numerical results show that the wells indeed affect each other in terms of performance due to pressure buildup. Depending on what Pareto optimal solution is selected by a decision maker, the agents will have different incentives to seek to form different coalitions. The stability of any coalition ultimately depends on how the value, i.e., the total amount of injected CO$_2$ is distributed within the coalition.

\section{Acknowledgements}
This work was financed in part by the Norwegian Research Council grant 336294,
Expansion of Resources for CO$_2$ Storage on the Horda Platform (ExpReCCS).	

\bibliographystyle{plain}
\bibliography{bibliography}

\end{document}